# Dependencies of prime numbers in a tuple

Victor Volfson

ABSTRACT. This paper investigates the dependence between primes in tuples through the analysis of the Hardy-Littlewood constant. A detailed analysis of the behavior of the constant for the pattern $(0, d)$ is conducted, depending on the arithmetic properties of $d$, including cases of convergence to the twin prime constant, divergence to infinity, and oscillatory behavior. Using methods from analytic number theory, the limiting distributions of the corresponding multiplicative functions are studied. It is shown that for symmetric tuples, the Hardy-Littlewood constant monotonically decreases as the tuple length decreases, indicating a weakening dependence between the primes. Theoretical conclusions are supported by calculations for specific symmetric tuples.

# 1. INTRODUCTION

A prime pattern (or template) of length $k$ is a set of $k$ distinct even numbers $\mathcal{H} = \{h_1, h_2, \ldots, h_k\}$, where is usually assumed to be $h_1 = 0$. This pattern defines the relative distances between prime numbers in a set of the form $\{n+h_1, n+h_2, \ldots, n+h_k\}$. Such a set is called a prime tuple.

Classic examples are the pattern $\{0, 2\}$, which corresponds to a tuple of twin primes (e.g. $11, 13$). The pattern $\{0, 2, 6, 8\}$, which corresponds to a prime quadruples (e.g., $5, 7, 11, 13$).

We exclude from consideration patterns $\mathcal{H}$ that form a complete residue system modulo for some prime number $p$. Such patterns can only correspond to a finite number of prime tuples (for example, the pattern $\{0, 1\}$ corresponds to a single tuple $\{2, 3\}$). The remaining patterns will be called admissible.

The Cramer model [1] is used for probabilistic modeling of the distribution of prime numbers. Its main assumption is that the events "the number $n$ is prime" for different are independent, and the probability of each such event is $1/\ln n$, which is motivated by the prime number theorem.

The probability that all numbers in a tuple $\{n+h_1, \ldots, n+h_k\}$ are prime within this model would be asymptotically equal to the product:

$$\prod_{i=1}^{k} \frac{1}{\ln(n+h_i)} \sim \frac{1}{(\ln n)^k}. \tag{1.1}$$

Therefore, the expected number of such $k$-tuples up to $x$ as $x \to \infty$ is asymptotically equal to:

$$\pi_k(x) \sim \frac{x}{(\ln x)^k}. \tag{1.2}$$

However, Cramer's model is a significant simplification and does not take into account two important features of the distribution of prime numbers:

1. Parity. All prime numbers are odd except 2. All numbers will be odd for a tuple with even differences $n+h_i$. The probability that a random odd number is prime is approximately $2/\ln n$. Taking this factor into account increases the expected number of $k$-tuples to:

$$\pi_k(x) \sim \frac{x}{2} \cdot \left(\frac{2}{\ln x}\right)^k = 2^{k-1} \cdot \frac{x}{(\ln x)^k}. \tag{1.3}$$

2. The influence of small primes. In reality, primes are not independent; their distribution exhibits complex correlations due to divisibility by small primes. This effect is taken into account by the first Hardy–Littlewood conjecture [2]. According to it, the number of prime $k$-tuples with an admissible pattern $\mathcal{H}$ up to $x$ is asymptotically equal to:

$$\pi_k(x, \mathcal{H}) \sim \mathfrak{S}(\mathcal{H}) \frac{x}{(\ln x)^k}, \tag{1.4}$$

where $\mathfrak{S}(\mathcal{H})$ is the Hardy–Littlewood constant, which depends on the pattern.

Comparing (1.2) and (1.4), we see that Cramer's model would assume $\mathfrak{S}(\mathcal{H}) = 1$. The actual constant is calculated using the formula:

$$\mathfrak{S}(\mathcal{H}) = \prod_p \left(1 - \frac{\nu_p(\mathcal{H})}{p}\right)\left(1 - \frac{1}{p}\right)^{-k}. \tag{1.5}$$

where the product is taken over all prime numbers $p$, and $\nu_p(\mathcal{H})$ is the number of distinct remainders modulo $p$ in the pattern $\mathcal{H}$. Thus, the constant $\mathfrak{S}(\mathcal{H})$ quantitatively characterizes the relationship between the events of "being prime" for numbers in a single tuple.

Many modern studies [3–5] rely on the assumption that the Hardy–Littlewood conjecture is true. An important result obtained in this direction is Gallagher's theorem [6], which states that, under the condition that the Hardy–Littlewood conjecture is true, the number of admissible prime $k$-tuples asymptotically obeys a Poisson distribution.

Having in mind an analysis of the degree of dependence between primes in a tuple is of considerable interest. Specifically, conjectures have been put forward that this dependence decreases with increasing tuple diameter $d(\mathcal{H}) = \max(\mathcal{H}) - \min(\mathcal{H})$, as well as, for a fixed diameter, with decreasing prime density in the interval under consideration.

Assuming the validity of the first Hardy–Littlewood conjecture, we pose the following tasks within the framework of this paper:

1. Investigate the behavior of the Hardy–Littlewood constant as a measure of the dependence of primes in a tuple.

2. Analyze conjectures that the dependence between primes in a tuple weakens with increasing diameter.

3. Investigate the nature of the dependence between primes in a tuple as the prime density decreases (for example, when considering longer intervals) for tuples of a fixed diameter.

The following chapters of this paper, assuming that the first Hardy–Littlewood conjecture is true, consider these questions.

2. CHANGE IN DEPENDENCE OF PRIMES IN A TUPLE WITH BAN INCREASE IN THE DIAMETR OF THE TUPLE

Let's consider the case where the length of a prime tuple is $k = 2$. The pattern of this tuple is $(0, d)$, where $d$ is the diameter of the tuple.

The pattern $H = (0, d)$ is always admissible for any $d$. This means that the remainders from dividing the pattern elements $0, d$ by $p$ do not cover the entire set of residues for any prime number $\{0, 1, \ldots, p - 1\}$.

For $p = 2$ elements: $0$ and $d$ have remainder $0 (\bmod 2) = 0$.

Since $d$ is even, the remainder is $d (\bmod 2) = 0$, therefore the number of different remainders is 1 and does not cover $\mathbb{Z}_2 = \{0, 1\}$.

For $p > 2$:

- if $p | d$, then the remainders are equal to $0$ and $0$, that is, the number of distinct

residues is equal to 1, which does not cover $\mathbb{Z}_p$;

- if $p \neq d$, then the remainders are equal to $0$ and $d (\bmod p) \neq 0$, that is, there are two distinct remainders, but $p \geq 3$, so they do not cover all $p$ residues.

It seems, that the dependence between primes in the tuple should decrease for $d \to \infty$ and therefore $C(H) \to 1$. However, this is not the case. The value $C(H) \to \infty$ for some sequences $d(n)$ at $n \to \infty$. The value $C(H)$ does not change for other sequences $d(n)$ at $n \to \infty$. Let us demonstrate this.

Based on (1.5), the Hardy-Littlewood constant for a tuple with pattern $H = (0,d)$ is calculated using the formula for pairs of prime numbers ($k = 2$):

$$C(H) = 2 \cdot \prod_{p>2} \frac{1 - \frac{v_H(p)}{p}}{(1 - \frac{1}{p})^2}, \qquad (2.1)$$

where $v_H(p)$ is the number of distinct remainders of the tuple elements modulo $p$.

Prime factors (except 2) for a tuple with the pattern $H = (0,d)$: $p = 3$, $p = 5, \ldots$, $p \leq d$.

If $p \mid d$, the remainders are the same ($v_H(p) = 1$). The remainders are different ($v_H(p) = 2$) for the rest $p > 2$. Therefore, the formula for the constant (2.1) simplifies to:

$$C(H) = 2C_2 \cdot \prod_{\substack{p>2 \\ p \mid d}} \frac{p-1}{p-2}, \qquad (2.2)$$

where $C_2 \approx 0.66016181584686957$ is the prime factor constant.

Using (2.2), we will review the limit of the Hardy-Littlewood constant for patterns of the form $H = (0, d(n))$ for various sequences $d(n)$ as $n \to \infty$.

Let us consider sequences with bounded prime divisors.

If all $d(n)$ have only a fixed set of prime divisors, then the product in formula (2.2) for $C(H)$ is bounded, and the limit is finite.

Example: $d(n) = 2^n$. Prime factors: only 2 (counted in the factor 2). Then the following holds $C(H) = 2C_2 \cdot 1 = 2C_2 \approx 1.32032$ and:

$$\lim_{n \to \infty} C(H) = 2C_2 \approx 1.32032. \qquad (2.3)$$

Now we examine sequences with prime divisors that grow sufficiently quickly.

If the prime factors grow so rapidly that the series $\sum \frac{1}{p}$ converges, then the product in (2.2) converges to a finite limit.

Example: $d(n) = p_n$ ($n$-th prime number). There is only one prime divisor $p_n$ (other than 2) for each $d(n)$. The factor in formula (2.2):

$$\frac{p_n - 1}{p_n - 2} \to 1$$

for $n \to \infty$.

Therefore, the following holds:

$$\lim_{n \to \infty} C(H) \to 2C_2 \approx 1.32032. \tag{2.4}$$

Now let's talk about a sequence with an increasing number of small prime factors.

If $d(n)$ includes more and more small primes, then the product in (2.2) diverges, and $C(H) \to \infty$.

Example: $d(n) = q_n\#$, where $q_n\#$ is the primorial (the product of the first prime numbers). We will use the asymptotics:

$$\prod_{\substack{p>2 \\ p \leq q}} \frac{p-1}{p-2} \sim 2C_2 \ln q_n \quad \text{(The proof is given below – assertion 1).}$$

Then, based on (2.2):

$$C(H) \sim 2C_2 \cdot (2C_2 \ln q) = 4C_2^2 \ln q_n. \tag{2.5}$$

Therefore, based on (2.5), we obtain:

$$\lim_{n \to \infty} C(H) = \infty. \tag{2.6}$$

Therefore, the dependence between primes in a tuple increases without bound in this case.

Let us study sequences with moderate growth of prime divisors.

If $d(n)$ has prime factors growing as $n$, then $C(H)$ can fluctuate or grow slowly.

Example: $d(n) = n$.

The value $C(H)$ depends on the number and size of prime factors $n$. If $n$ has many small prime factors, $C(H)$ can be large. If $n$ is a prime number, then the value can be little. There is no single limit; the behavior is unstable.

Thus, based on (2.3), (2.4), (2.6), the constant $C(H)$ for the pattern $H=(0,d(n))$ can converge to $2C_2$, diverge to infinity, or fluctuate depending on the arithmetic properties of $d(n)$.

The dependence between primes in the corresponding pattern $H=(0,d(n))$ tuples behaves similarly when $d(n)\to\infty$.

Assertion 1

$$L(q)=\prod_{\substack{p>2\\p\leq q}}\frac{p-1}{p-2}\sim 2C_2\ln q, \qquad (2.7)$$

where $C_2$ is the constant of twin primes.

Proof

Let's consider the logarithm:

$$\ln\left(1+\frac{1}{p-2}\right)=\frac{1}{p-2}-\frac{1}{2(p-2)^2}+\frac{1}{3(p-2)^3}-\cdots$$

We use the expansion of the logarithm into a series:

$$L(q)=\ln\ln q+K+o(1).$$

Then:

$$L(q)=\sum_{\substack{p>2\\p\leq q}}\frac{1}{p-2}-\frac{1}{2}\sum_{\substack{p>2\\p\leq q}}\frac{1}{(p-2)^2}+\frac{1}{3}\sum_{\substack{p>2\\p\leq q}}\frac{1}{(p-2)^3}-\cdots=\frac{1}{p-2}-\frac{1}{2(p-2)^2}+\frac{1}{3(p-2)^3}-\cdots$$

Let us denote:

$$S_1(q)=\sum_{\substack{p>2\\p\leq q}}\frac{1}{p-2},\quad S_2(q)=\sum_{\substack{p>2\\p\leq q}}\frac{1}{(p-2)^2},\quad S_3(q)=\sum_{\substack{p>2\\p\leq q}}\frac{1}{(p-2)^3},\quad\ldots$$

Note that:

$$S_1(q) = \sum_{\substack{p>2 \\ p \leq q}} \frac{1}{p} + \sum_{\substack{p>2 \\ p \leq q}} \left( \frac{1}{p-2} - \frac{1}{p} \right).$$

It is known that:

$$\sum_{p \leq q} \frac{1}{p} = \ln \ln q + B + o(1),$$

where $B \approx 0.261497$ is the Meissel-Mertens constant.

The second series converges:

$$\sum_{p>2} \left( \frac{1}{p-2} - \frac{1}{p} \right) = \sum_{p>2} \frac{2}{p(p-2)} = L,$$

where $L$ is some constant.

Thus:

$$S_1(q) = \ln \ln q + B + L + o(1).$$

The series $S_k(q)$ converge to constants as $k \geq 2$:

$$S_2(q) = M_2 + o(1), \quad S_3(q) = M_3 + o(1), \quad \ldots$$

where $M_k = \sum_{p>2} \frac{1}{(p-2)^k}$.

Hence:

$$L(q) = \ln \ln q + B + L - \frac{1}{2} M_2 + \frac{1}{3} M_3 - \cdots + o(1).$$

Let us denote the constant:

$$K = B + L - \frac{1}{2} M_2 + \frac{1}{3} M_3 - \cdots.$$

Then we get:

$$L(q) = \ln \ln q + K + o(1),$$

and

$$\prod_{\substack{p>2 \\ p\leq q}} \frac{p-1}{p-2} = e^{L(q)} \sim e^K \ln q.$$

It remains to show that $e^K = 2C_2$.

This follows from a well-known result on the asymptotics of this product, where the constant $2C_2$ arises from the definition of the constant of prime twins:

$$C_2 = \prod_{p>2}\left(1 - \frac{1}{(p-1)^2}\right).$$

Thus, the following holds:

$$\prod_{\substack{p>2 \\ p\leq q}} \frac{p-1}{p-2} \sim 2C_2 \ln q,$$

which corresponds to (2.7).

Now let us consider the function $f(d) = C(0,d)$ defined for even values $d$, where $C(0,d)$ is the Hardy-Littlewood constant for the pattern of primes of the form $(0,d)$.

Based on (2.2), this function is determined by the formula:

$$f(d) = 2C_2 \cdot \prod_{\substack{p>2 \\ p\mid d}} \frac{p-1}{p-2}, \tag{2.8}$$

where $C_2$ is the constant of twin primes.

Having in mind (2.8), we form a normalized function:

$$g(d) = \frac{f(d)}{2C_2} = \prod_{\substack{p>2 \\ p\mid d}} \frac{p-1}{p-2} \tag{2.9}$$

for even $d \geq 2$.

Let us extend this function to all natural numbers for the convenience of studying this function and define an arithmetic function $h(n)$ for $n \geq 1$:

$$h(n) = \prod_{\substack{p\mid n \\ p>2}} \frac{p-1}{p-2}, \tag{2.10}$$

while we further define $h(1) = 1$ and if $p = 2$ for $k \geq 1$, then $h(2^k) = 1$, and if $p > 2$, then $h(p^k) = \dfrac{p-1}{p-2}$.

The function $h(n)$ is multiplicative, since for coprime numbers $h(nm) = h(n)h(m)$.

The following theorem of Wintner is known [7].

Let $h(n)$ is a multiplicative arithmetic function and $b(n) = h * \mu(n)$, where $*$ is the Dirichlet convolution.

Then, if

$$\sum_{n=1}^{\infty} \frac{|b(n)|}{n} < \infty,$$

then the mean value:

$$M(h) = \lim_{x \to \infty} \frac{1}{x} \sum_{n \leq x} h(n)$$

exists and is equal to:

$$M(h) = \sum_{n=1}^{\infty} \frac{b(n)}{n}. \qquad (2.11)$$

Based on (2.10), in our case:

$$h(n) = \prod_{\substack{p|n \\ p>2}} \frac{p-1}{p-2}.$$

Values are $b(1) = 1$, $b(p) = h(p) - 1 = \dfrac{p-1}{p-2} - 1 = \dfrac{1}{p-2}$, $b(p^k) = h(p^k) - h(p^{k-1}) = 0$ for $k \geq 2$.

Based on (2.11) and the above, the average value $h(n)$ is:

$$M(h) = \sum_{n=1}^{\infty} \frac{b(n)}{n} = \prod_{p} (1 + \frac{b(p)}{p}) = \prod_{p>2} (1 + \frac{1}{p(p-2)}) = \frac{1}{C_2} = 1{,}51478... \qquad (2.12)$$

Thus, $h(n)$ is not limited, but it is limited on average (2.12).

Now, let's define the moment of the $k$-th order of the arithmetic function $h(n)$:

$$M_k(h) = \lim_{x \to \infty} \frac{1}{x} \sum_{n \leq x}^{\infty} h^k(m). \tag{2.13}$$

We will use the fact that $h^k(m)$ is also a multiplicative arithmetic function for which the Halberstam-Richter theorem holds [8].

Let $H(n) \geq 0$ is a multiplicative arithmetic function for which the condition holds $\sum_p \frac{H(p)-1}{p} < \infty$, then it is true:

$$M(H) = \prod_p (1-\frac{1}{p})(1+\frac{H(p)}{p}+\frac{H(p^2)}{p^2}+...). \tag{2.14}$$

It is true in our case:

$$H(n) = h^k(n) = \prod_{\substack{p>2 \\ p|k}} (\frac{p-1}{p-2})^k. \tag{2.15}$$

Let's define - $H(2^n) = 1$.

Let's check the conditions of the Halberstam-Richter theorem. $H(n) \geq 0$ is a multiplicative function, $H(p)-1 = (\frac{p-1}{p-2})^k - 1$.

Let us expand (2.15) in a series for large $p$:

$$\frac{p-1}{p-2} = 1 + \frac{1}{p-2} = 1 + \frac{1}{p} + O(\frac{1}{p^2}),$$

$$(\frac{p-1}{p-2})^k = (1+\frac{1}{p}+O(\frac{1}{p^2}))^k = 1 + \frac{k}{p} + O(\frac{1}{p^2}).$$

Therefore:

$$\sum_p \frac{H(p)-1}{p} \sim \sum_p \frac{k}{p^2} < \infty. \tag{2.16}$$

Thus, all conditions of the Halberstem-Richter theorem are fulfilled and, based on (2.14), the average value $H(n)$ can be determined:

$$M(\mathrm{H}) = \prod_p (1-\frac{1}{p})(1+\frac{H(p)}{p}+\frac{H(p^2)}{p^2}+...).$$

Contribution of value $H(2^n)=1$ for $p=2$ to the formula:

$$(1-\frac{1}{2})(1+\frac{1}{2}+\frac{1}{2^2}+...) = \frac{1}{2} \cdot 2 = 1.$$

Let's substitute the value $H(p^n) = (\frac{p-1}{p-2})^k$ for $p \geq 3$ into the formula:

$$1+\frac{H(p)}{p}+\frac{H(p^2)}{p^2}+... = 1+(\frac{p-1}{p-2})^k(\frac{1}{p}+\frac{1}{p^2}+...) = 1+(\frac{p-1}{p-2})^k \frac{1}{p-1}.$$

Now we multiply by $(1-\frac{1}{p}) = \frac{p-1}{p}$ and get:

$$(1-\frac{1}{p})(1+\frac{H(p)}{p}+\frac{H(p^2)}{p^2}+...) = \frac{p-1}{p}(1+(\frac{p-1}{p-2})^k \frac{1}{p-1}) = \frac{p-1}{p}+\frac{1}{p}(\frac{p-1}{p-2}).$$

Therefore, the mean value $H(n)$ is:

$$M(\mathrm{H}) = \prod_p (\frac{p-1}{p}+\frac{1}{p}(\frac{p-1}{p-2})^k) = \prod_p (1-\frac{1}{p}+\frac{h^k(p)}{p}). \qquad (2.17)$$

The meaning is $h^k(2)=1$ for $p=2$.

The meaning (for $p>2$) is:

$$h^k(p) = (\frac{p-1}{p-2})^k \approx 1+\frac{k}{p^k}. \qquad (2.18)$$

Based on (2.17) and (2.18) we obtain:

$$M_2[f] = (2C_2)^2 M_2[h] = 4C_2^2 M_2[h] = (1,320)^2 \, 2,649 = 4,619. \qquad (2.19)$$

Based on (2.19), the variance is:

$$D(h) = M_2(h) - M_1^2(h) \approx 0,354. \qquad (2.20)$$

We explore $h(n)$ for extreme values.

The maximum value $h(n)$ is achieved on the segment $[1, x]$, as mentioned above, with values $n$ equal to the primorial:

$$\max_{n \leq x} h(n) = \prod_{2 < p \leq P} \frac{p-1}{p-2} = \prod_{2 < p \leq P} (1 + \frac{1}{p-2}). \qquad (2.21)$$

We take the logarithm of (2.21):

$$\ln(\prod_{2 < p \leq P}(1+\frac{1}{p-2})) = \sum_{2 < p \leq P} \ln(1+\frac{1}{p-2}) \sim \sum_{2 < p \leq P} \frac{1}{p-2} = \sum_{2 < p \leq P} \frac{1}{p} + O(1) = \ln \ln P + O(1). \qquad (2.22)$$

Since $n = \prod_{p \leq P} p \leq x$, then by the prime number theorem $n = \prod_{p \leq P} p \sim e^P \leq x$ or $P \sim \ln x$, therefore $\ln \ln P \sim \ln \ln \ln x$, then based on (2.22), the asymptotics of the maximum value is equal to:

$$\max_{n \leq x} h(n) \sim e^{\ln \ln \ln x + C} = e^C \ln \ln x. \qquad (2.23)$$

The value of the constant in (2.23) is $C = \ln 2$, therefore:

$$\max_{n \leq x} h(n) \sim 2 \ln \ln x. \qquad (2.24)$$

The minimum value $h(n) = 1$ is achieved at $n = 2^k$.

Now let's get back to the function $f(d) = 2C_2 h(d/2)$.

The mean value of the function is:

$$M_1[f] = 2C_2 M_1[h] = 2C_2 \frac{1}{C_2} = 2. \qquad (2.25)$$

Second-order moment:

$$M_2[f] = (2C_2)^2 M_2[h] = 4C_2^2 M_2[h] = (1,320)^2 2,649 = 4,619. \qquad (2.26)$$

Based on (2.25) and (2.26), the variance is:

$$D[f] = M_2[f] - M_1^2[f] \approx 4,619 - 4 = 0,619. \qquad (2.27)$$

Now let's consider the limiting distribution for the multiplicative arithmetic function $h(n)$. To do this, let's consider the limiting distribution of the additive arithmetic function $\ln(h(n))$.

The classical Erdős–Wintner theorem states [8] that a real additive arithmetic function $f(n)$ has a limit distribution if and only if the following three series over all primes $p$ converge simultaneously:

$$\sum_{|f(p)|>1} \frac{1}{p}, \sum_{|f(p)|<1} \frac{f(p)}{p}, \sum_{|f(p)|<1} \frac{f^2(p)}{p}. \tag{2.28}$$

Let us consider the additive arithmetic function $f(n) = \ln(h(n))$.

Since the logarithm of a product is equal to the sum of the logarithms, we obtain:

$$f(n) = \ln(h(n)) = \sum_{\substack{p \mid n \\ p>2}} \ln\left(\frac{p-1}{p-2}\right). \tag{2.29}$$

Based on (2.29), the asymptotics of this function on prime numbers has the form:

$$f(p) = \ln\left(\frac{p-1}{p-2}\right) = \ln(1 + \frac{1}{p-2}) \sim \frac{1}{p}. \tag{2.30}$$

Therefore:

$$f^2(p) \sim \frac{1}{p^2}. \tag{2.31}$$

Having in mind the obtained asymptotics (2.30), (2.31), we will check the conditions of the Erdős-Wintner theorem.

Since the condition $f(p) > 1$ is not satisfied for any prime $p$, the first series $\sum_{|f(p)|>1} \frac{1}{p}$ does not have a single term and converges trivially.

It follows from the asymptotics $f(p) \sim \frac{1}{p}$ that the general term of the second series behaves as $\frac{1}{p^2}$. It is known that the series $\sum_{p} \frac{1}{p^2}$ converges. Therefore, the original series also converges absolutely.

It follows from the asymptotics $f^2(p) \sim \dfrac{1}{p^2}$ that the general term of the third series behaves as $\dfrac{1}{p^2} \cdot \dfrac{1}{p} = \dfrac{1}{p^3}$. The series $\sum_p \dfrac{1}{p^3}$ - converges. Therefore, the third series also converges.

Thus, all three conditions of the Erdős-Wintner theorem are satisfied. This proves that the additive function $f(n) = \ln(h(n))$ has a limit distribution.

The distribution function, according to this theorem, will be an infinite product over all prime numbers:

$$F(t) = \prod_p \left(1 - \dfrac{1}{p}\right)\left(1 + \sum_{k=1}^{\infty} p^{-k} \exp\left(it f^k(p)\right)\right). \qquad (2.32)$$

According to the Erdős-Wintner theorem, the distribution of an additive arithmetic function is discrete if and only if the series $\sum_{f(p)\neq 0} \dfrac{1}{p}$ converges. Otherwise, if the series diverges, the distribution is continuous (and even absolutely continuous).

Let's consider the series:

$$\sum_{f(p)\neq 0} \dfrac{1}{p} = \sum_{p>2} \dfrac{1}{p}. \qquad (2.33)$$

It is known that the series $\sum_p \dfrac{1}{p}$ over all primes diverges (this follows from the divergence of the harmonic series and the asymptotics $\sum_{p \leq x} \dfrac{1}{p} \sim \ln \ln x$). Since the series over all primes diverges, the series over odd primes $\sum_{p>2} \dfrac{1}{p}$ also diverges, since it differs from the complete series by a finite term of 1/2.

Since the series $\sum_{f(p)\neq 0} \dfrac{1}{p}$ diverges, then, in accordance with the Erdős-Wintner theorem, the distribution of the function $f(n) = \ln(h(n))$ is absolutely continuous.

A multiplicative arithmetic function:

$$h(n) = \prod_{\substack{p \mid n \\ p>2}} \dfrac{p-1}{p-2}$$

has a discrete distribution over the set of natural numbers. This follows from its multiplicative nature and the fact that it takes values from a specific discrete set.

The distribution $h(n)$ is discrete and asymmetric (skewed to the right).

Most values $h(n)$ are close to 1 (for example, for numbers without odd prime factors or with large prime factors).

Rare values $h(n)$ can be large (when n contains many small odd prime factors), which increases the mean.

The median is less than the mathematical expectation because the distribution is skewed to the right.

The function $f(d) = 2C_2 h(d/2)$ also has a discrete distribution similar to $h(n)$, but with the different means and variances discussed above.

3. RELATIONSHIPS BETWEEN PRIMES IN SYMMETRIC TUPLES

Definition: A symmetric prime tuple of length $k$ with diameter $d$ is a prime tuple that has a pattern $H = (h_1, h_2, \ldots, h_k)$ with the condition $h_i + h_{k-i} = d$.

This means that the elements of a tuple consisting of prime numbers are symmetrically distributed about the center.

Assertion 2

Let $\mathcal{H} = \{h_1, h_2, \ldots, h_k\}$ is the pattern of a symmetric tuple of diameter $d$. Then, the following holds for the Hardy-Littlewood constants and any of its proper subpatterns $H'$:

$$\mathfrak{S}(\mathcal{H}') < \mathfrak{S}(\mathcal{H}). \tag{3.1}$$

Proof

Hardy-Littlewood constant for pattern $\mathcal{H}$ length $k$ is:

$$\mathfrak{S}(\mathcal{H}) = \prod_p \left(1 - \frac{\nu_p(\mathcal{H})}{p}\right)\left(1 - \frac{1}{p}\right)^{-k}. \tag{3.2}$$

Based on (3.2), the Hardy-Littlewood constant for a subpattern $\mathcal{H}$ of length $k' < k$:

$$\mathfrak{S}(\mathcal{H}') = \prod_p \left(1 - \frac{\nu_p(\mathcal{H}')}{p}\right)\left(1 - \frac{1}{p}\right)^{-k'}. \tag{3.3}$$

Let's choose $p_0 > d$. Then for:

- $p \geq p_0$: $v_p(\mathcal{H}) = k, v_p(\mathcal{H}') = k'$, since all elements are distinct in absolute value $p$;

- $p < p_0$: the values of $v_p(\mathcal{H})$ and $v_p(\mathcal{H}')$ depend on the structure of the tuples.

Based on (3.2) and (3.3) we obtain the relation:

$$\frac{\mathfrak{S}(\mathcal{H}')}{\mathfrak{S}(\mathcal{H})} = \underbrace{\prod_{p<p_0} \frac{1-\frac{v_p(\mathcal{H}')}{p}}{1-\frac{v_p(\mathcal{H})}{p}} \left(1-\frac{1}{p}\right)^{k-k'}}_{P_{\text{small}}} \cdot \underbrace{\prod_{p \geq p_0} \frac{1-\frac{k'}{p}}{1-\frac{k}{p}} \left(1-\frac{1}{p}\right)^{k-k'}}_{P_{\text{big}}}. \qquad (3.4)$$

The multiplier in (3.4) for $p \geq p_0$ has the form:

$$A_p = \frac{p-k'}{p-k} \cdot \left(1-\frac{1}{p}\right)^{k-k'}. \qquad (3.5)$$

Based on (3.5), the following holds in the case $k' = k-1$ (when deleting one element):

$$A_p = \frac{p-(k-1)}{p-k} \cdot \left(1-\frac{1}{p}\right) = 1 + \frac{k-1}{p(p-k)}.$$

Therefore, the product $P_{\text{big}} = \prod_{p \geq p_0} A_p$ converges to a finite number $1+\delta$ where $\delta > 0$.

Now let's consider the product at $p < p_0$:

$$P_{\text{small}} = \prod_{p<p_0} B_p, \qquad B_p = \frac{p-v_p(\mathcal{H}')}{p-v_p(\mathcal{H})} \cdot \left(1-\frac{1}{p}\right). \qquad (3.6)$$

Now we will prove the lemma, and then we will complete the proof of the statement using it.

Lemma

It is true for a symmetric pattern $\mathcal{H}$ and any of its subpatterns $\mathcal{H}'$:

$$P_{\text{small}} \leq \theta < 1, \qquad (3.7)$$

where the constant $\theta$ depends only on $k$ and $p_0$, but not on the specific choice of $\mathcal{H}'$.

Proof

Since $p_0$ is chosen as $p_0 > d$, where $d$ is the diameter of the tuple for the pattern $\mathcal{H}$, the set of primes $p < p_0$ is finite. This means that the product $P_{small}$ includes a finite number of factors.

The value $v_p(\mathcal{H})$ (the number of distinct residues modulo $p$) for each prime $p < p_0$ is determined by the structure of the symmetric tuple pattern $\mathcal{H}$.

Symmetry imposes strict restrictions on the possible values of $v_p(\mathcal{H})$, for example, it is always true $v_2(\mathcal{H}') = 1$ for $p = 2$.

Values $v_p(\mathcal{H}')$ can be either equal to $v_p(\mathcal{H})$ or less than for a subpattern $\mathcal{H}'$.

In particular:

- it is always true - $v_2(\mathcal{H}') = 1$ for $p = 2$, since all elements of the tuple are odd. Therefore, $B_2 = \dfrac{2-1}{2-1} \cdot \left(1 - \dfrac{1}{2}\right) = \dfrac{1}{2}$.

- it is hold $v_p(\mathcal{H}') = v_p(\mathcal{H})$ for other small primes $p$, for subtuples $\mathcal{H}'$, what gives $B_p = 1 \cdot \left(1 - \dfrac{1}{p}\right) < 1$.

Let's consider the maximum $P_{small}$ over all possible subpatterns $\mathcal{H}'$. Since the number of subpatterns is finite, the maximum is achieved.

Let us show that this maximum is strictly less than 1:

- the product $P_{small}$ always contains the factor $B_2 = \dfrac{1}{2}$;

- the factors $B_p$ are bounded above for other primes $p < p_0$. For example, if $v_p(\mathcal{H}') = v_p(\mathcal{H}) - 1$, then

$$B_p = \dfrac{p - v_p(\mathcal{H}) + 1}{p - v_p(\mathcal{H})} \cdot \left(1 - \dfrac{1}{p}\right) = \left(1 + \dfrac{1}{p - v_p(\mathcal{H})}\right) \cdot \left(1 - \dfrac{1}{p}\right).$$

Since $v_p(\mathcal{H}) \leq p - 1$ (from the admissibility of the tuple), then $p - v_p(\mathcal{H}) \geq 1$, and

$$B_p \leq 2 \cdot \left(1 - \frac{1}{p}\right) = 2 - \frac{2}{p}.$$

This gives $B_3 \leq \frac{4}{3}$ for $p = 3$, $B_3 \leq 1,2$ for $p = 5$, etc.

In general, for product $P_{small}$:

$$P_{small} \leq \frac{1}{2} \prod_{3 \leq p < p_0} \left(2 - \frac{2}{p}\right) = \prod_{3 \leq p < p_0} (1 - \frac{1}{p}) < 1.$$

Therefore, it is hold $P_{small} \leq \theta < 1$ for any $p_0$.

Thus, the lemma is proved: there exists a constant $\theta < 1$, depending only on $k$ and $p_0$, such that $P_{small} \leq \theta$ for all $\mathcal{H}' \subset \mathcal{H}$.

Now we complete the proof of assertion 2.

We have:

$$\frac{\mathfrak{S}(\mathcal{H}')}{\mathfrak{S}(\mathcal{H})} = P_{small} \cdot P_{big} \leq \theta \cdot (1 + \delta).$$

We choose $p_0$ large enough so that $\theta \cdot (1 + \delta) < 1$. This is always possible based on (3.7), since:

- $\theta$ is bounded by a constant independent of $p_0$;

- $\delta \to 0$ for $p_0 \to \infty$.

Hence:

$$\frac{\mathfrak{S}(\mathcal{H}')}{\mathfrak{S}(\mathcal{H})} < 1 \quad \Rightarrow \quad \mathfrak{S}(\mathcal{H}') < \mathfrak{S}(\mathcal{H}).$$

Corollary 1 (Monotonicity under symmetric reduction of a tuple)

Let $\mathcal{H}_k$ is the pattern of a symmetric tuple of length $k$. If $\mathcal{H}_{k-2}$ is obtained from $\mathcal{H}_k$ by deleting a symmetric pair of elements $\{h, h'\}$ (where $h$ and $h'$ are symmetric about the center), then:

$$\mathfrak{S}(\mathcal{H}_{k-2}) < \mathfrak{S}(\mathcal{H}_k).$$

Proof

The pattern $\mathcal{H}_{k-2}$ is a subpattern of the symmetric tuple $\mathcal{H}_k$. By Assertion 2, the strict inequality holds for any subpattern of the symmetric tuple. Since $\mathcal{H}_{k-2} \subset \mathcal{H}_k$, then $\mathfrak{S}(\mathcal{H}_{k-2}) < \mathfrak{S}(\mathcal{H}_k)$.

This consequence allows us to "jump" over even lengths while maintaining symmetry and monotonicity.

Corollary 2 (Monotonicity under successive symmetric deletion)

Suppose we have a sequence of tuples

$$\mathcal{H}_m \subset \mathcal{H}_{m+1} \subset \cdots \subset \mathcal{H}_k,$$

where each tuple $\mathcal{H}_i$ is symmetric, and each subsequent tuple is obtained from the previous one by adding a symmetric pair of elements. Then:

$$\mathfrak{S}(\mathcal{H}_m) < \mathfrak{S}(\mathcal{H}_{m+1}) < \cdots < \mathfrak{S}(\mathcal{H}_k).$$

Proof

Let's consider any two adjacent subpatterns of a symmetric tuple: $\mathcal{H}_i \subset \mathcal{H}_{i+1}$.

By construction:

- $\mathcal{H}_{i+1}$ is a subpattern of a symmetric tuple;

- $\mathcal{H}_i$ is its subpattern (obtained by removing the symmetric pair).

Based on assertion 2: $\mathfrak{S}(\mathcal{H}_i) < \mathfrak{S}(\mathcal{H}_{i+1})$.

Since this is true for every pair of adjacent tuples in the chain, we obtain monotonicity:

$$\mathfrak{S}(\mathcal{H}_m) < \mathfrak{S}(\mathcal{H}_{m+1}) < \cdots < \mathfrak{S}(\mathcal{H}_k).$$

Note: The wording says "adding a symmetrical pair", which is equivalent to "removing a symmetrical pair" when moving in the opposite direction.

Conclusions:

1. Corollary 1 allows us to compare the Hardy-Littlewood constants for symmetric tuples of different even lengths.

2. Corollary 2 provides a powerful tool for analyzing the asymptotic behavior of the constants as the size of symmetric tuples increases.

3. Both corollaries preserve the structural property of symmetry, which is important for many applications in number theory (e.g., in conjectures on the density of primes in tuples).

Thus, the Hardy-Littlewood constant for symmetric tuples decreases monotonically as $k$.

Based on the data obtained, the table below was compiled, showing the change in the Hardy-Littlewood constant $C(H)$ as the tuple length decreases for the pattern $k = 17, d = 240$.

Table

| Tuple length ($k$) | Pattern | $C(H)$ | Relationship $C(H_{k-2})/C(H_k)$ |
|---|---|---|---|
| 17 | [0, 6, 24, 36, 66, 84, 90, 114, 120, 126, 150, 156, 174, 204, 216, 234, 240] | $2.0427 \times 10^8$ | - |
| 15 | [0, 24, 36, 66, 84, 90, 114, 120, 126, 150, 156, 174, 204, 216, 240] | $1.5265 \times 10^7$ | 0.0747 |
| 13 | [0, 36, 66, 84, 90, 114, 120, 126, 150, 156, 174, 204, 240] | $1.6893 \times 10^6$ | 0.1107 |
| 11 | [0, 66, 84, 90, 114, 120, 126, 150, 156, 174, 240] | $1.1665 \times 10^5$ | 0.0691 |
| 9 | [0, 84, 90, 114, 120, 126, 150, 156, 240] | $9.6778 \times 10^3$ | 0.0830 |

| 7 | [0, 90, 114, 120, 126, 150, 240] | $5.0374 \times 10^2$ | 0.0521 |
| 5 | [0, 114, 120, 126, 240] | $9.2634 \times 10^1$ | 0.1839 |
| 3 | [0, 120, 240] | $1.1433 \times 10^1$ | 0.1234 |
| 2 | [0, 240] | $3.5209 \times 10^0$ | 0.3079 |

This table clearly demonstrates that as the number of elements in the tuple decreases, the Hardy-Littlewood constant $C(H)$ decreases monotonically, which is consistent with theoretical expectations about the weakening of the dependence between primes in symmetric tuples.

## 4. CONCLUSION AND SUGGESTIONS FOR FURTHER WORK

Next article will continue to study the asymptotic behavior of some arithmetic functions.

## 5. ACKNOWLEDGEMENTS

Thanks to everyone who has contributed to the discussion of the paper.